\documentclass[a4paper,10pt]{amsart}
\usepackage{color}

\usepackage[utf8]{inputenc} 
\usepackage[T1]{fontenc} 
\usepackage[francais,british]{babel} 
\usepackage{graphicx} 
\usepackage{csquotes}
\usepackage{amsmath,amsfonts,amssymb,amsthm} 
\usepackage{listings}
\usepackage{nccmath}
\usepackage{subfigure}
\usepackage{enumerate}
\usepackage{mathrsfs}
\usepackage[backend=bibtex,isbn=false,url=false,doi=false,style=authoryear-comp,firstinits=true,maxbibnames=9]{biblatex}
\usepackage[margin=1.5in]{geometry}
\newtheorem{theorem}{Theorem}[section]

\newtheorem{property}[theorem]{Property}
\newtheorem{fact}[theorem]{Fact}

\theoremstyle{definition}
\newtheorem{ex}[theorem]{Example}

\newtheorem{defn}[theorem]{Definition}

\theoremstyle{remark}
\newtheorem{remark}[theorem]{Remark}

\newcommand{\p} {\ensuremath {\mathbb{P}}}
\newcommand{\E} {\ensuremath {\mathbb{E}}}

\newcommand{\N} {\ensuremath {\mathbb{N}}}
\newcommand{\R} {\ensuremath {\mathbb{R}}}

\newcommand{\Z} {\ensuremath {\mathbb{Z}}}
\newcommand{\I} {\ensuremath {\mathbb{I}}}
\newcommand{\F} {\ensuremath {\mathscr{F}}}

\newcommand{\Nn} {\ensuremath {\mathscr{N}}}

\newcommand{\A} {\ensuremath {\mathscr{A}}}
\newcommand{\X} {\ensuremath {\mathscr{X}}}
\newcommand{\Ti} {\ensuremath {\mathscr{T}}}

\newcommand{\B} {\ensuremath {\mathscr{B}}}
\newcommand{\M} {\ensuremath {\mathscr{M}}}

\newcommand{\mo} {\ensuremath {\mathscr{P}}}
\newcommand{\Qi} {\ensuremath {\mathscr{Q}}}

\newcommand{\Wh} {\ensuremath {\mathscr{W}}}

\title{Le Cam theory on the comparison of statistical models}
\author{Ester ~Mariucci}
\thanks{The research leading to these results has received funding from the European Research Council 
under ERC Grant Agreement 320637.}
\address{Leiden University.}
\email{ester.mariucci@gmail.com}
\date{\today}
\begin{document}
\begin{abstract}
We recall the main concepts of the Le Cam theory of statistical experiments, especially the notion of Le Cam distance and its properties. We also review classical tools for bounding such a distance before presenting some examples. A proof of the classical equivalence result between density estimation problems and Gaussian white noise models will be analyzed.
\end{abstract}

\maketitle


\section{Introduction}
The theory of \emph{Mathematical Statistics} is based on the notion of \emph{statistical model}, also called \emph{statistical experiment} or just \emph{experiment}. A statistical model, as in its original formulation due to \textcite{B51}, is a triple
$$\mo=(\Omega,\Ti,(P_\theta: \theta\in \Theta)),$$
where $(\Omega,\Ti)$ is a sample space, $\Theta$ is a set called the \emph{parameter space} and $(P_\theta :\theta\in \Theta)$ is a family of probability measures on $(\Omega,\Ti)$. This definition is a mathematical abstraction intended to represent a concrete experiment; consider for example the following situation taken from the book of \textcite{LC2000}.
A physicist decides to estimate the half life of Carbon $14$, $C^{14}$. He supposes that the life of a $C^{14}$ atom has an exponential distribution with parameter $\theta$ and, in order to develop his investigation, he takes a sample of $n$ atoms of $C^{14}$. The physicist fixes in advance the duration of the experiment, say $2$ hours, and then he counts the number of disintegrations. Formally, this leads to the definition of the statistical model $\mo_1=(\N,\mathcal{P}(\N), (P_\theta:\theta\in(0,\infty)))$ where $P_\theta$ represents the law of the random variable $X$ counting the number of disintegrations observed in $2$ hours. This is not the only way to proceed if we want to estimate the half life of Carbon $14$. Indeed, the physicist could choose to consider the first random time $Y$ after which a fixed number of disintegrations, say $10^6$, have occurred. In this case he will represent the experiment via the statistical model $\mo_2=(\R_+,\B(\R_+), (Q_\theta:\theta\in(0,\infty)))$ where $Q_\theta$ 
is the law of the random variable $Y$. A natural question is then how much ``statistical information'' the considered experiments contain or, more precisely, when the experiment $\mo_1$ will be more informative than $\mo_2$ and conversely. 


The quest for comparison of statistical experiments was initiated by the paper of \textcite{BSS} followed by the papers of \textcite{B51,B53} where the following definition was introduced: ``$\mo_1$ is more informative than $\mo_2$'' if for any bounded loss function $L$, $\|L\|_{\infty}\leq 1$, and any decision procedure $\rho_2$ in the experiment $\mo_2$ there exists a decision procedure $\rho_1$ in the experiment $\mo_1$ such that 
$$R_\theta(\mo_1,\rho_1,L)\leq R_\theta(\mo_2,\rho_2,L),\quad \forall \theta\in\Theta.$$
Here we denote by $R_\theta(\mo_1,\rho_1,L)$ and $R_\theta(\mo_2,\rho_2,L)$ the \emph{statistical risk} for the experiments $\mo_1$ and $\mo_2$, respectively.

However, this can lead to two models being non-comparable. This issue was solved by Le Cam  who introduced the notion of deficiency $\delta(\mo_1,\mo_2)$. We will give a precise definition in the forthcoming sections. Here, we only remark two interesting properties:
\begin{itemize}
 \item $\delta(\mo_1,\mo_2)$ is a well defined non-negative real number for every two given statistical models $\mo_1$ and $\mo_2$ sharing the same parameter space.
 \item For every loss function $L$ with $0\leq L\leq 1$ and every decision procedure $\rho_2$ available on $\Theta$ using $\mo_2$, there exists a decision procedure $\rho_1$ in $\mo_1$ such that for all $\theta\in\Theta$,
 $$R_\theta(\mo_1,\rho_1,L)\leq R_\theta(\mo_2,\rho_2,L)+\delta(\mo_1,\mo_2).$$
\end{itemize}
This solves the issue mentioned above: It could be that both $\delta(\mo_1,\mo_2)$ and $\delta(\mo_2,\mo_1)$ are strictly positive, in which case they will not be comparable according to the first definition; nevertheless, we can still say ``how much information'' we lose when passing from one model to the other one.
Le Cam's theory has found applications in several problem in statistical decision theory and it has been developed, for example, for nonparametric regression, nonparametric density estimation problems, generalized linear models, diffusion models, Lévy models, spectral density estimation problem. 
Historically, the first results of asymptotic equivalence in a nonparametric context date from 1996 and are due to \textcite{BL} and \textcite{N96}.
The first two authors have shown the asymptotic equivalence of nonparametric regression and a Gaussian white noise model while the third one those of density estimation problems and Gaussian white noise models. Over the years many generalizations of these results have been proposed such as \textcite{BL2002,GN2002,ro04,C2007,cregression,R2008,C2009,R2013,schmidt14} for nonparametric regression or \textcite{cmultinomial,j03,BC04,esterdensity} for nonparametric density estimation models.
Another very active field of study is that of diffusion experiments. The first result of equivalence between diffusion models and Euler scheme was established in 1998, see \textcite{NM}. In later papers generalizations of this result have been considered (see \textcite{C14, esterdiffusion}) as well as different statistical problems always linked with diffusion processes (see, e.g., \textcite{D,CLN,R2006,rmultidimensionale}).
Among others we can also cite equivalence results for generalized linear models (see, e.g., \textcite{GN}), time series (see, e.g., \textcite{GN2006,NM}), GARCH model (see, e.g., \textcite{B}), functional linear regression (see, e.g.,\textcite{M2011}), spectral density estimation (see, e.g. \textcite{GN2010}), volatility estimation (see, e.g. \textcite{R11}) and jump models (see, e.g.,  \textcite{esterESAIM,esterlevy}). Negative results are somewhat harder to come by; the most notable among them are \textcite{sam96,B98,wang02}.
Another new research direction that has been explored involves quantum statistical experiments (see, e.g., \textcite{BF}).

The aim of this survey paper is to present some basic concepts of the Le Cam theory of asymptotic equivalences between statistical models.
Our aim in this review is to give an accessible introduction to the subject. Therefore, we will not follow the most general approach to the theory, also because such an approach is already available in the literature, see e.g., \textcite{lecam,LC2000,vdv2002}.
In order to achieve such a goal, the paper has been organized as follows. In Section \ref{sec:lecam} we recall the definition of the Le Cam 
distance and its statistical meaning. Particular attention has been payed to the interpretation of the Le Cam distance in terms of decision 
theory. In Section \ref{sec:bounds} we collect some classical tools to control the Le Cam distance before passing to some examples described in 
Section \ref{sec:examples}. Section \ref{sec:density} is devoted to show in details a proof of a classical result in Le Cam theory, namely the 
asymptotic equivalence between density estimation problems and Gaussian white noise models.  

\section{Deficiency and Le Cam distance}\label{sec:lecam}
As we have already pointed out, a possible way to compare two given statistical models (having the same parameter space) could be to compare the corresponding risk functions or to ask ``how much information'' we lose when passing from one model to the other one, saying that there is no loss if we have at our disposal a mechanism able to convert the observations from the distribution $P_{1,\theta}$ to observations from $P_{2,\theta}$. If we adopt the latter point of view a natural formalization for such a mechanism is the notion of Markov kernel.

\begin{defn}
 Let $(\X_i,\Ti_i)$, $i=1,2$, be two measurable spaces. A \emph{Markov kernel} $K$ with source $(\X_1,\Ti_1)$ and target $(\X_2,\Ti_2)$ is a map 
 $K: \X_1\times\Ti_2\to [0,1]$ with the following properties:
 \begin{itemize}
  \item The map $x\mapsto K(x, A)$ is $\Ti_1$-measurable for every $A\in\Ti_2$.
  \item The map $A\mapsto K(x,A)$ is a probability measure on $(\X_2,\Ti_2)$ for every $x\in \X_1$.
 \end{itemize}
\end{defn}

We will denote by $K:(\X_1,\Ti_1)\to(\X_2,\Ti_2)$ a Markov kernel with source $(\X_1,\Ti_1)$ and target $(\X_2,\Ti_2)$.
 
Starting from a Markov kernel $K:(\X_1,\Ti_1)\to(\X_2,\Ti_2)$ and a probability measure $P_1$ on $(\X_1,\Ti_1)$ one can construct a probability measure on $(\X_2,\Ti_2)$ in the following way:
$$K P_1(A)=\int K(x,A)P_1(dx),\quad \forall A\in\Ti_2.$$

Roughly speaking we can think that two models $\mo_1$ and $\mo_2$ contain ``the same amount of information about $\theta$'' if there exist two Markov
kernels, $K_1$ and $K_2$, not depending on $\theta$, such that $K_1P_{1,\theta}=P_{2,\theta}$ and $K_2P_{2,\theta}=P_{1,\theta}$. This idea has been 
formalized in the sixties by Lucien Le Cam and led to the notion of the deficiency, hence to the introduction of a pseudo-metric on the class of all statistical experiments having the same parameter space.

The definition of the deficiency in its most general form involves the notion of ``transition'' which is a generalization of the concept of Markov 
kernel. In this paper, however, we prefer to keep things simpler and only focus on the case in which one has to deal with dominated statistical models
having Polish sample spaces (see below for a definition). The advantage is that in this case the definition of deficiency simplifies and the abstract concept of transition coincides with that of Markov kernel (see Proposition 9.2 in \textcite{N96}).

\begin{defn}
 A statistical model $\mo_1=(\X_1,\Ti_1,(P_{1,\theta}:\theta\in\Theta))$ is called \emph{Polish} if its sample space $(\X_1,\Ti_1)$ is a separable 
 completely metrizable topological space. 
 
 $\mo_1$ is said to be \emph{dominated} if there exists a $\sigma$-finite measure $\mu$ on $(\X_1,\Ti_1)$ such that, for all $\theta\in\Theta$, $P_{1,\theta}$
 is absolutely continuous with respect to $\mu$. The measure $\mu$ is called the \emph{dominating measure}. 
\end{defn}

\begin{ex}
 Typical examples of Polish spaces in probability theory are the spaces $\R,\R^n,\R^{\infty}$, the space $C_T$ of continuous functions on $[0,T]$
 equipped with the supremum norm $d(x,y)=\sup_{0\leq t\leq T}|x_t-y_t|$, the space $D$ of càdlàg functions equipped with the Skorokhod metric.
\end{ex}

\begin{defn}
 Let $Q_1$ and $Q_2$ be two probability measures defined on a measurable space $\Omega$. The \emph{total variation distance} between $Q_1$ and $Q_2$
 is defined as the quantity:
 $$\|Q_1-Q_2\|_{TV}=\sup_{A\subseteq \Omega}|Q_1(A)-Q_2(A)|=\frac{1}{2}L_1(Q_1,Q_2),$$
where $L_1(Q_1,Q_2)$ denotes the $L_1$ norm between $Q_1$ and $Q_2$.
 \end{defn}

 \begin{defn}
 Let $\mo_i=(\X_i,\Ti_i,(P_{i,\theta}:\theta\in\Theta))$, $i=1,2$, be two experiments.
 The \emph{deficiency} $\delta(\mo_1,\mo_2)$ of $\mo_1$ with respect to $\mo_2$ is the number
  $$\delta(\mo_1,\mo_2)=\inf_T\sup_{\theta\in\Theta}\|TP_{1,\theta}-P_{2,\theta}\|_{TV},$$
  for an infimum taken over all Markov kernels $T: (\X_1,\Ti_1)\to(\X_2,\Ti_2)$ and $\|\cdot\|_{TV}$ denotes the total variation distance.
  \end{defn}
  \begin{defn}
  The \emph{Le Cam distance} or \emph{$\Delta$-distance} between $\mo_1$ and $\mo_2$ is defined as
  $$\Delta(\mo_1,\mo_2)=\max(\delta(\mo_1,\mo_2),\delta(\mo_2,\mo_1)).$$
  \end{defn}
  
  The $\Delta$-distance is a pseudo-metric on the space of all statistical models: It satisfies the triangle inequality 
  $\Delta(\mo_1,\mo_3)\leq\Delta(\mo_1,\mo_2)+\Delta(\mo_2,\mo_3)$ but the equality $\Delta(\mo_1,\mo_2)=0$ does not imply that 
  $\mo_1$ and $\mo_2$ actually coincide. 
  
Concerning the glossary, when $\delta(\mo_1,\mo_2)=0$ (i.e. if the experiment $\mo_2$ can be reconstructed from the experiment $\mo_1$ by a Markov kernel), we will say that
$\mo_2$ is \emph{less informative} than $\mo_1$, or that $\mo_1$ is \emph{better} than $\mo_2$, or that $\mo_1$ is \emph{more informative} than $\mo_2$.
When $\Delta(\mo_1,\mo_2)=0$ the models $\mo_1$ and $\mo_2$ are said to be \emph{equivalent} and two sequences of statistical models $(\mo_{1,n})_{n\in\N}$ 
and $(\mo_{2,n})_{n\in \N}$ are called \emph{asymptotically equivalent} when $\Delta(\mo_{1,n},\mo_{2,n})\to 0$ as $n\to\infty$.

A way to interpret the Le Cam distance between experiments is to see it as a numerical indicator of the cost needed to 
reconstruct one model from the other one and vice-versa, via Markov kernels.
But, as we said in the introduction, a way to compare statistical models that seems just as natural is to compare the respective risk functions. Let us then highlight how the definition of the deficiency has a clear interpretation in terms of statistical decision theory.
To that aim, we will start by recalling the standard framework:

\begin{itemize}
 \item A \emph{statistical model}, which is just an indexed set $\mo=(\X,\Ti,(P_\theta:\theta\in \Theta))$ of probability measures all defined on the same measurable space $(\X,\Ti)$, for some set $\X$ equipped with a $\sigma$-field $\Ti$. The elements of $\Theta$ are sometimes called the \emph{states of Nature}.
 \item A space $A$ of possible actions or decisions that the statistician can take after observing $x\in\X$. For example, in estimation problems we can take $A=\Theta$. To make sense of the notion of integral on $A$ we need it to be equipped with a $\sigma$-field $\mathcal{A}$.
 \item A loss function $L:\Theta\times A\mapsto (-\infty,\infty]$, with the interpretation that action $z\in A$ incurs a loss $L(\theta,z)$ when $\theta$ is the true state of Nature. 
   \item A \emph{(randomized) decision rule} $\rho$ in $\mo$ is a Markov kernel $\rho:(\X,\Ti)\to(A,\mathscr A)$.
  \item The \emph{risk} is:
  $$R_\theta(\mo,\rho,L)=\int_{\X}\bigg(\int_A L(z,\theta)\rho(y,dz)\bigg)P_\theta(dy).$$
More precisely, the standard interpretation of risk is as follows.
The statistician observes a value $x\in\X$ obtained from a probability measure $P_{\theta}$. He does not know the value of $\theta$ and must take a decision  $z\in A$. He does so by choosing a probability measure $\rho(x,\cdot)$ on $A$ and picking a point in $A$ at random according to $\rho(x,\cdot)$. If he has chosen $z$ when the true distribution of $x$ is $P_{\theta}$, he suffers a loss $L(\theta,z)$. His average loss when $x$ is observed is then $\int L(\theta,z)\rho(x,dz)$. His all over average loss when $x$ is picked according to $P_{\theta}$ is the integral $\int\big(\int L(\theta,z)\rho(x,dz)\big)P_{\theta}(dx)$.

  \end{itemize}

A very important result allowing to translate the notion of deficiency as described above in a decision theory language is the following:
\begin{theorem}[See \textcite{lecamsufficiency} or Theorem 2, page 20 in \textcite{lecam}]\label{teo:lecam}
Let $\varepsilon>0$ be fixed. $\delta(\mo_1,\mo_2)< \varepsilon$ if and only if: $\forall$ decision rule $\rho_2$ on $\mo_2$ and for all 
bounded loss function $L$, $\|L\|_\infty\leq 1$, there exists a decision rule $\rho_1$ on $\mo_1$ such that
$$R_\theta(\mo_1,\rho_1,L)< R_\theta(\mo_2,\rho_2,L)+\varepsilon,\quad \forall \theta\in\Theta.$$
\end{theorem}
In other words we have that
\begin{align*}
\delta(\mo_1,\mo_2) &=\inf_{\rho_1}\sup_{\rho_2}\sup_{\theta}\sup_{L}|R(\mo_1,\rho_1,L,\theta)-R(\mo_2,\rho_2,L,\theta)|,
\end{align*}
where the last supremum is taken on the set of all loss functions $L$ s.t. $0\leq L(\theta,z)\leq 1$, $\forall z\in A$, $\forall \theta\in\Theta$
and $\rho_i$ belongs to the set of all randomised decision procedures in the experiment $\mo_i$, $i=1,2$.
\begin{remark}
 An important consequence of the previous theorem is that if two sequences of experiments $(\mo_{1,n})_{n\in\mathbb N}$ and 
 $(\mo_{2,n})_{n\in\mathbb N}$ are asymptotically equivalent in the Le Cam's sense then asymptotic properties of any inference problem are the same 
 for these experiments. 
This means that when two sequences of statistical experiments are proven to be asymptotically equivalent it is enough to choose the simplest 
one, to study there the inference problems one is interested in and to transfer the knowledge about such inference problems to the more complicated 
sequence, via Markov kernels.
\end{remark}
\subsection{How to transfer decision rules via randomisations}

\textcolor{white}{p}

Let $\mo_{i,n}=(\X_{i,n},\mathscr T_{i,n}, (P_{i,n,\theta}:\theta\in\Theta))$, $i=1,2$, be two sequences of statistical models sharing the
same parameter space $\Theta$ and having Polish sample spaces $(\X_{i,n},\mathscr T_{i,n})$. 
Suppose that there exist Markov kernels $K_n$ such that $\|K_nP_{1,n,\theta}-P_{2,n,\theta}\|_{TV}\to 0$ uniformly on the parameter space. 
Then, given a decision rule (or an estimator) $\pi_{2,n}$ on $\mo_{2,n}$ we can define a decision rule $\pi_{1,n}$ on $\mo_{1,n}$ that, 
asymptotically, has the same statistical risk as $\pi_{2,n}$. To show that let us start by considering the easier case in which both $K_n$ and $\pi_{2,n}$ are deterministic.    
More precisely, we suppose that $K_n$ is of the form $K_n(A)=\I_{A}S_n(x)$ for all $A\in \mathscr{T}_{2,n}$ for some functions $S_n$. 

Then, we have (suppressing the index $n$ to shorten notations): 

\begin{align*}
 &\bigg|\int_{\X_1}L(\theta,\pi_1(y))P_{1,\theta}(dy)-\int_{\X_2}L(\theta,\pi_{2}(y))P_{2,\theta}(dy)\bigg|\leq \\
&\bigg|\int_{\X_{1}}L(\theta,\pi_{1}(y))P_{1,\theta}(dy)-\int_{\X_{2}}L(\theta,\pi_{2}(y))KP_{1,\theta}(dy)\bigg|\\
& +\bigg|\int_{\X_{2}}L(\theta,\pi_{2}(y))\big[KP_{1,\theta}(dy)-P_{2,\theta}(dy)\big]\bigg| \\
&\leq \bigg|\int_{\X_1}L(\theta,\pi_1(y))P_{1,\theta}(dy)-\int_{\X_1}L(\theta,\pi_{2}(S(y)))P_{1,\theta}(dy)\bigg|+\|L\|_{\infty}\|KP_{1}-P_2\|_{TV}
 \end{align*}
 In particular, assuming that the loss function $L$ is bounded by $1$ and defining 
 $$\pi_1(y):=\pi_2(S(y))$$
 one finds that
 $$\bigg|\int_{\X_1}L(\theta,\pi_1(y))P_{1,\theta}(dy)-\int_{\X_2}L(\theta,\pi_{2}(y))P_{2,\theta}(dy)\bigg|\leq \|KP_{1}-P_2\|_{TV}\to 0,$$
 that is, the decision rule $\pi_{1,n}(y)=\pi_{2,n}(S_n(y))$ has asymptotically the same risk as $\pi_{2,n}$. The same kind of computations 
 work in the general case in which the $K_n$'s are not deterministic and $(\pi_{2,n})$ is a sequence of decision rule having $(A_n,\A_n)$ as action's spaces. In this case one can show that the randomized sequence of decision rules
 $$\pi_{1,n}(y,C):=\int_{\X_{2,n}}\pi_{2,n}(x,C)K(y,dx),\quad \forall y\in\X_{1,n},\ \forall C\in\A_n$$
 has asymptotically the same risk as $\pi_{2,n}$.



\begin{remark}\label{independentkernels}
Let $P_i$ be a probability measure on $(E_i,\mathcal{E}_i)$ and $K_i$ a Markov kernel on $(G_i,\mathcal G_i)$. One can then define a Markov kernel $K$ on $(\prod_{i=1}^n E_i,\otimes_{i=1}^n \mathcal{G}_i)$ in the following way:
 $$K(x_1,\dots,x_n; A_1\times\dots\times A_n):=\prod_{i=1}^nK_i(x_i,A_i),\quad \forall x_i\in E_i,\ \forall A_i\in \mathcal{G}_i.$$
 Clearly $K\otimes_{i=1}^nP_i=\otimes_{i=1}^nK_iP_i$.
\end{remark}

\section{How to control the Le Cam distance}\label{sec:bounds}
Even if the definition of deficiency has a perfectly reasonable statistical meaning, it is not easy to compute: Explicit computations have appeared but they are rare (see \textcite{torgersen72,torgersen74,hansen74} and Section 1.9 in \textcite{spok}).
More generally, one may hope to find more easily some upper bounds for the $\Delta$-distance. We collect below some useful techniques for this purpose.

\begin{property}\label{delta0}
 Let $\mo_j=(\X,\Ti,(P_{j,\theta}; \theta\in\Theta))$, $j=1,2$, be two statistical models having the same sample space and define 
 $\Delta_0(\mo_1,\mo_2):=\sup_{\theta\in\Theta}\|P_{1,\theta}-P_{2,\theta}\|_{TV}.$
 Then, $\Delta(\mo_1,\mo_2)\leq \Delta_0(\mo_1,\mo_2)$.
\end{property}
In particular, Property \ref{delta0} allows us to bound the $\Delta$-distance between statistical models sharing the same sample space by means of classical bounds for the total variation distance. To that aim, we collect below some useful (and classical) results.

\begin{fact}[see \cite{LC69}, p. 35]
 Let $P_1$ and $P_2$ be two probability measures on $\X$, dominated by a common measure $\xi$, with densities $g_{i}=\frac{dP_{i}}{d\xi}$, $i=1,2$. Define
 \begin{align*}
  L_1(P_1,P_2)&=\int_{\X} |g_{1}(x)-g_{2}(x)|\xi(dx), \\
  H(P_1,P_2)&=\bigg(\int_{\X} \Big(\sqrt{g_{1}(x)}-\sqrt{g_{2}(x)}\Big)^2\xi(dx)\bigg)^{1/2}.
 \end{align*}
Then,
\begin{equation*} 
 \frac{H^2(P_1,P_2)}{2}\leq \|P_1-P_2\|_{TV}=\frac{1}{2}L_1(P_1,P_2)\leq H(P_1,P_2).
\end{equation*}
\end{fact}
An important property is the following:

\begin{property}\label{hellprodotto}
If $\mu$ and $\nu$ are product measures defined on the same measurable space, $\mu=\bigotimes_{j=1}^{m}\mu_{j}$ and $\nu=\bigotimes_{j=1}^{m}\nu_{j}$, then
\begin{equation*}
H^2(\mu,\nu)=2\bigg[1-\prod_{j=1}^{m}\bigg[1-\frac{H^2(\mu_j,\nu_j)}{2}\bigg]\bigg].
\end{equation*}
\end{property}
\begin{proof}
See, e.g., \textcite{zolotarev83}, p. 279. 
\end{proof}

Thus one can express the distance between distributions of vectors with independent components in terms of the component-wise distances. A consequence of Property \ref{hellprodotto} is:
\begin{property}\label{h}
If $\mu$ and $\nu$ are product measures defined on the same measurable space, $\mu=\bigotimes_{j=1}^{m}\mu_{j}$ and $\nu=\bigotimes_{j=1}^{m}\nu_{j}$, then
 $$H^2(\mu,\nu)\leq \sum_{i=1}^m H^2(\mu_i,\nu_i).$$
\end{property}
\begin{proof}
See, e.g., \textcite{strasser}, Lemma 2.19. 
\end{proof}

\begin{property}\label{hellnormale}
The Hellinger distance between two normal distributions $\mu\sim \mathcal{N}(m_1,\sigma_1^2)$ and $\nu\sim \mathcal{N}(m_2,\sigma_2^2)$ is:
\begin{align*}
H^2(\mu,\nu)&=2\bigg[1-\bigg[\frac{2\sigma_1\sigma_2}{\sigma_1^2+\sigma_2^2}\bigg]^{1/2}\exp\bigg[-\frac{(m_1-m_2)^2}{4(\sigma_1^2+\sigma_2^2)}\bigg]\bigg]\\
 &\leq 2\bigg|1-\frac{\sigma_1^2}{\sigma_2^2}\bigg|+\frac{(m_1-m_2)^2}{2\sigma_2^2}.
\end{align*}
  \end{property}
\begin{proof}
 See, e.g., \textcite{esterESAIM}, Fact 1.5.
\end{proof}

\subsection{The likelihood process}
Another way to control the Le Cam distance lies in the deep relation linking the equivalence between experiments to the proximity of the distributions of the related likelihood ratios.

Let $\mo_j=(\X_j,\Ti_j,(P_{j,\theta}:\theta\in\Theta))$ be a statistical model dominated by $P_{j,\theta_0}$, $\theta_0\in \Theta$, and let $\Lambda_j(\theta)=\frac{dP_{j,\theta}}{dP_{j,\theta_0}}$ be the density of $P_{j,\theta}$ with respect to $P_{j,\theta_0}$. In particular, one can see  $\Lambda_j(\theta)$ as a real random variable defined on the probability space $(\X_j,\Ti_j)$, i.e. one can see $(\Lambda_j(\theta))_{\theta\in\Theta}$ as a stochastic process. For that reason we introduce the notation $\Lambda_{\mo_j}:=(\Lambda_j(\theta),\theta\in\Theta)$ and we call $\Lambda_{\mo_j}$ the \emph{likelihood process}.

A key result of the theory of Le Cam is the following.
\begin{property}
Let $\mo_j=(\X_j,\Ti_j,(P_{j,\theta}:\theta\in\Theta))$, $j=1,2$, be two experiments. If the family $(P_{j,\theta}:\theta\in\Theta)$ is dominated by $P_{j,\theta_0}$, then $\mo_1$ and $\mo_2$ are equivalent if and only if their likelihood processes under the dominating measures $P_{1,\theta_0}$ and $P_{2,\theta_0}$ coincide.
\end{property}
\begin{proof}
see \textcite{strasser}, Corollary 25.9. 
\end{proof}
Let us now suppose that there are two processes $(\Lambda_j^{n,*}(\theta))_{\theta\in \Theta}$, $j=1,2$ defined on a same probability space $(\X^*,\Ti^*,\Pi^*)$ and such that the law of $(\Lambda_j^{n}(\theta))_{\theta\in \Theta}$ under $P_{j,\theta_0}$ is equal to the law of $(\Lambda_j^{n,*}(\theta))_{\theta\in \Theta}$ under $\Pi^*$, $j=1,2$. Then, the following holds (see \textcite{LC2000}, Lemma 6).

\begin{property}
 If $\Lambda_{\mo_1}$ and $\Lambda_{\mo_2}$ are the likelihood processes associated with the experiments $\mo_1$ and $\mo_2$, then
 $$\Delta(\mo_1^n,\mo_2^n)\leq \sup_{\theta\in \Theta}\E_{\Pi^*}\Big|\Lambda_1^{n,*}(\theta)-\Lambda_2^{n,*}(\theta)\Big|.$$
\end{property}

\subsection{Sufficiency and Le Cam distance}
A very useful tool, when comparing statistical models having different sample spaces, is to look for a sufficient statistic.   
The introduction of the term \emph{sufficient statistic} is usually attributed to R.A. Fisher who gave several definitions of the concept. We cite here the presentation of the subject from \textcite{lecamsufficiency}. Fisher's most relevant statement seems to be the requirement ``...that the statistic chosen should summarize the whole of the relevant information supplied by the sample.'' Such a requirement may be made precise in various ways, the following three interpretations are the most common.
\begin{enumerate}[(i)]
 \item \emph{The classical, or operational definition of sufficiency}, claims that a statistic $S$ is sufficient if, given the value of $S$, one can proceed to a post-experimental randomization reproducing variables which have the same distributions as the originally observable variables.
 \item \emph{The Bayesian interpretation}. A statistic $S$ is sufficient if for every a priori distribution of the parameter the a posteriori distributions of the parameter given $S$ is the same as if the entire result of the experiment was given.
 \item \emph{The decision theoretical concept.} A statistic $S$ is sufficient if for every decision problem and every decision procedure made available by the experiment there is a decision procedure, depending on $S$ only, which has the same performance characteristics.
\end{enumerate}
The study of sufficiency in an abstract way can be found in \textcite{HS}. The last section of such a work is named ``The value of sufficient statistics in statistical methodology" and starts with the following observation:
\blockquote{
\emph{We gather from conversations with some able and prominent mathematical statisticians that there is doubt and disagreement about just what a sufficient statistic is sufficient to do, and in particular about in what sense if any it contains ``all the information in a sample"}.}

In  \textcite{B54} a continuation of the work of \textcite{HS} can be found. A particular effort was done to highlight the interest of using sufficient statistics in statistical methodology. One of the main results in \textcite{B54} is Theorem 7.1 establishing the equivalence of the decision theoretical concept of sufficiency and the operational concept in terms of conditional probabilities. We mention this fact because of its similarity with the result of Le Cam, here stated as Theorem \ref{teo:lecam}, that is the core of the theory of comparison of statistical experiments.

Formally, let $\mo=(\X,\Ti,(P_{\theta}: \theta\in\Theta))$ be a statistical model. A \emph{statistic} is a measurable map from a measurable space $(\X,\Ti)$ to another measurable one $(\X_2,\Ti_2)$. We denote by $S_{\#}P_\theta$ the image law of $S$ under $P_\theta$, i.e
$S_{\#}P_\theta(B)=P_\theta(S^{-1}(B))$, for all $B\in \Ti_2$.
 \begin{defn}
 $S$ is a \emph{sufficient statistic} for $(P_{\theta}: \theta\in\Theta)$ if for any $A\in \Ti$ there exists a function $\phi_A$, with $\phi_A \circ S$ $\Ti$-measurable, such that
 $$P_\theta(A\cap S^{-1}(B))=\int_B\phi_A(y)S_{\#}P_\theta(dy), \quad \forall A\in \Ti,\ \forall B\in \Ti_2,\ \forall \theta\in\Theta.$$
 An arbitrary subalgebra $\Ti_0$ of $\Ti$ is said to be \emph{sufficient} for $(P_\theta:\theta\in\Theta)$ if for all $A\in \Ti$ there exists a $\Ti_0$-measurable function $\phi_A$ such that
 $$P_\theta(A\cap A_0)=\int_{A_0}\phi_A(x)P_\theta(dx),\quad \forall A_0\in\Ti_0,\ \forall \theta\in\Theta.$$
 The set $\{S^{-1}(B): B\in \Ti_2\}$ is called the \emph{subalgebra induced by the statistic $S$}.
\end{defn}
\begin{property}[See, e.g. \textcite{B54}]
A statistic $S$ is sufficient for $(P_\theta:\theta\in \Theta)$ if the subalgebra induced by $S$ is sufficient for $(P_\theta:\theta\in \Theta)$.
\end{property}
In accordance with the notation introduced in Section \ref{sec:lecam}, we will state Theorem 7.1 in \textcite{B54} as follows (recall that $(A,\A)$ denotes the \emph{action/decision space}.)
\begin{theorem}[See Theorem 7.1, \textcite{B54}]
If the subalgebra $\Ti_0$ of $\Ti$ is sufficient for $(P_\theta:\theta\in \Theta)$, then for every decision rule $\rho:(\X,\Ti)\mapsto (A,\A)$ there exists a decision rule $\pi:(\X,\Ti_0)\mapsto (A,\A)$ 
such that
$$P_\theta\rho(C)=P_\theta\pi(C), \quad \forall C\in \A,\ \forall \theta\in \Theta.$$
\end{theorem}

Before focusing on the relation between the notion of sufficient statistic and the one of equivalence between statistical models,
let us recall the Neyman-Fisher factorization theorem, a powerful tool for identifying sufficient statistics for a given dominated family 
of probabilities.
Let $(P_\theta:\theta\in\Theta)$ be a family of probabilities on $(\Omega,\Ti)$, absolutely continuous with respect to a $\sigma$-finite measure $\mu$, 
and denote by $p_\theta:=\frac{dP_\theta}{d\mu}$ the density.
 \begin{theorem}
  A statistic $S:(\Omega,\Ti)\to(\X,\B)$ is sufficient for $(P_\theta:\theta\in\Theta)$ if and only if there exists a $\B$-measurable function $g_\theta$ $\forall \theta\in\Theta$ and a $\Ti$-measurable function $h\neq 0$ such that
  $$p_\theta(x)=g_\theta(S(x))h(x),\quad \mu\text{-a.s. } \forall x\in\Omega.$$
 \end{theorem}
An important result linking the Le Cam distance with the existence of a sufficient statistic is the following: 
 \begin{property}\label{fatto3}
 Let $\mo_i=(\X_i,\Ti_i,(P_{i,\theta}: \theta\in\Theta))$, $i=1,2$, be two statistical models. 
Let $S:\X_1\to\X_2$ be a sufficient statistic
such that the distribution of $S$ under $P_{1,\theta}$ is equal to $P_{2,\theta}$. Then $\Delta(\mo_1,\mo_2)=0$. 
\end{property}

\begin{proof}
 In order to prove that $\delta(\mo_1,\mo_2)=0$ it is enough to consider the Markov kernel $M:(\X_1,\Ti_1)\to(\X_2,\Ti_2)$ defined as 
 $M(x,B):=\I_B(S(x))$ $\forall x\in\X_1$ and $\forall B\in\Ti_2$. Conversely, to show that $\delta(\mo_2,\mo_1)=0$ one can consider the Markov kernel 
 $K:(\X_2,\Ti_2)\to (\X_1,\Ti_1)$ defined as 
 $K(y,A)=\E_{P_{2,\theta}}(\I_A|S=y)$, $\forall y\in\X_2$, $\forall A\in\Ti_1.$
 Since $S$ is a sufficient statistic, the Markov kernel $K$ does not depend on $\theta$. Denoting by $S_{\#}P_{1,\theta}$ the distribution of $S$ under $P_{1,\theta}$, one has:
 $$KP_{2,\theta}(A)=\int K(y,A)P_{2,\theta}(dy)=\int \E_{P_{2,\theta}}(\I_A|S=y)S_{\#}P_{1,\theta}(dy)=P_{1,\theta}(A).$$
 \end{proof}

For asymptotic arguments, one also needs an appropriate version of the notion of sufficiency.
\begin{defn}
Let $\mo_n=(\X_n,\Ti_n,(P_{n,\theta}: \theta\in\Theta))$ be a sequence of statistical models. The sequence of subalgebras $\tilde \Ti_n$ of $\Ti_n$ is \emph{asymptotically sufficient} for $(P_{n,\theta}:\theta\in\Theta)$ if $\Delta(\mo_n,\mo_n{_{|_{\tilde\Ti_n}}})\to 0$, where $\mo_n{_{|_{\tilde\Ti_n}}}$ denotes the restriction of the experiment $\mo_n$ to $\tilde \Ti_n$, i.e. $\mo_n{_{|_{\tilde\Ti_n}}}=(\X_n,\tilde \Ti_n, (\tilde P_{n,\theta}:\theta\in \Theta))$, where $\tilde P_{n,\theta}(A)=P_{n,\theta}(A),$ for all $A\in\tilde\Ti_n.$
\end{defn}
This is a stronger notion than asymptotic equivalence; indeed, let $\mo_{1,n}$ and $\mo_{2,n}$ be two sequences of experiments having the same parameter space. Then, by the triangle inequality, it is clear that if there exist two sequences $S_{1,n}$ and $S_{2,n}$ of asymptotically sufficient statistics in $\mo_{1,n}$ and $\mo_{2,n}$ respectively, taking values in the same measurable space, and such that  
$$\sup_{\theta\in \Theta}\|{S_{1,n}}_{\#}P_{1,\theta}-{S_{2,n}}_{\#}P_{2,\theta}\|_{TV}\to 0\quad \text{as }n\to\infty,$$
then the sequences $\mo_{1,n}$ and $\mo_{2,n}$ are asymptotically equivalent.
We also recall that an important generalization of the notion of sufficiency is the notion of insufficiency. The discussion of this concept is beyond the purposes of this paper, the reader is referred to \textcite{lecam74,} or Chapter 5 in \textcite{lecam} for an exhaustive treatment of the subject.

\section{Examples}\label{sec:examples}
To better understand what is the typical form of an asymptotic equivalence result let us analyze some examples. As a toy example let us start by considering the following parametric case. 
\begin{ex}
Let $\mo_{1,n}$ be the statistical model associated with the observation of a vector $X$ of $n$ independent Gaussian random variables $\Nn(\theta,1)$. Here the inference concerns $\theta $ and the parameter space $\Theta$ will be an interval of $\R$. Formally
$$\mo_{1,n}=(\R^n,\B(\R^n), (P_{1,\theta}:\theta\in\Theta)),$$
where $P_\theta$ is the law of $X$.

Then, let us denote by $\mo_{2,n}$ the experiment associated with the observation of the empirical mean relative to the previous random variables, i.e. 

$$\mo_{2,n}=(\R,\B(\R), (P_{2,\theta}:\theta\in\Theta)),$$
where $P_{2,\theta}$ is the law of a Gaussian random variable $\Nn(\theta,1/n)$.
By means of the Neyman-Fisher factorization theorem it is easy to see that the application $S:\R^n\to \R$, $S(x_1,\dots,x_n)=\frac{\sum_{i=1}^n x_i}{n}$ is a sufficient statistic. An immediate application of Property \ref{fatto3} implies that $\Delta(\mo_{1,n},\mo_{2,n})=0$ for all $n$. 
\end{ex}
Before passing to some examples in a nonparametric framework, let us recall a result due to Carter and concerning the asymptotic equivalence between a multinomial and a Gaussian multivariate experiment. The parameter space will be a subset of $\R^m$ and the reason for which we focus on such a result
lies on its being a very useful tool in establishing global asymptotic equivalence results for density estimation problems.

\begin{ex}\label{ex:carter}
 Let $X=(X_1,\dots, X_m)$ be a random vector having a multinomial distribution of parameters $n$ and $(p_1,\dots,p_m)$ with $p_i\geq 0$ for $i=1,
 \dots,m$ and $\sum_{i=1}^m p_i=1$. 

Denote by $\mo$ the statistical model associated with a multinomial distribution $P_\theta=\mathcal M(n;(\theta_1,\dots,\theta_m))$ with parameters 
$\theta=(\theta_1,\dots,\theta_m)$ that belong to $\Theta_R\subset \R^m$, a set consisting of all vectors of probabilities such that 
$$\frac{\max_i \theta_i}{\min_i \theta_i}\leq R.$$
The main result in \textcite{cmultinomial} is a bound of the Le Cam distance between statistical models associated with
multinomial distributions and multivariate normal distributions with the same means and covariances as the multinomial ones.
More precisely, let us denote by $\Qi$ the statistical model associated with a family of multivariate normal distributions $Q_\theta=\Nn(\mu, \Sigma)$, 
$\theta\in \Theta_R$, where 
$$\mu=(n\theta_1,\dots,n\theta_m), \Sigma=(\sigma_{i,j})_{i,j=1,\dots,m}\text{ with }
\sigma_{i,j}=n\theta_i(1-\theta_i)\delta_{i=j}-n\theta_i\theta_j\delta_{i\neq j}.$$

\begin{theorem}[see \textcite{cmultinomial}, p. 709]\label{cartermultinomial}
 With the notations above,
 $$\Delta(\mo,\Qi)\leq C_R\frac{m\ln m}{\sqrt n}$$
 for a constant $C_R$ that depends only on $R$.
\end{theorem}
Another interesting result contained in \textcite{cmultinomial} is the approximation of $\Qi$ by a Gaussian experiment with independent coordinates. Let 
us denote by $\tilde \Qi$ the statistical model associated with $m$ independent Gaussian random variables $\Nn(\sqrt{\theta_i},1/(4n))$, $i=1,\dots,m$. 
\begin{theorem}[see \textcite{cmultinomial}, p. 717--719]\label{cartergaussian}
 With the notations above,
 $$\Delta(\Qi,\tilde \Qi)\leq C_R\frac{m}{\sqrt n}$$
 for a constant $C_R$ that depends only on $R$.
\end{theorem}
\end{ex}
Let us now consider some examples in a nonparametric framework. More precisely, we will recall the results of \textcite{BL} and \textcite{N96} that are the first asymptotic equivalence results for nonparametric experiments.
\begin{ex}
In \textcite{BL}, the authors consider the problem of estimating the function $f$ from a continuously observed Gaussian process $y(t)$, $t\in[0,1]$, which satisfies the SDE
 $$dy_t=f(t)dt+\frac{\sigma(t)}{\sqrt n} dW_t,\quad t\in[0,1],$$
 where $dW_t$ is a Gaussian white noise.
They find that the statistical model associated with the continuous observation of $(y_t)$ is asymptotically equivalent to the statistical model associated with its discrete counterpart, i.e. the nonparametric regression:
$$y_i=f(t_i)+\sigma(t_i)\xi_i,\quad i=1,\dots,n.$$
The time grid is uniform, $t_i=\frac{i-1}{n}$, and the $\xi_i$'s are standard normal variables; they assume that $f$ varies in a nonparametric subset $\F$ of $L_2[0,1]$ defined by some smoothness conditions and $n$ tends to infinity not too slowly. More precisely, the drift function $f(\cdot)$ is unknown and such that, for $B$ a positive constant, one has:
$$\sup\big\{|f(t)|:t\in[0,1],f\in \F\big\}=B<\infty.$$
Moreover, defining
$$
\bar{f}_n(t)=\left\{
\begin{array}{ll}
f\big(\frac{i}{n}\big)&\textnormal{if} \quad \frac{i-1}{n}\leq t<\frac{i}{n},\quad i=1,\dots,n;\\
f(1) & \textnormal{if}\quad t=1;
\end{array}\right.
$$
one asks:
\begin{equation*}
\lim_{n\to\infty}\sup_{f\in\mathcal{F}}n\int_{0}^{1}\frac{(f(t)-\bar{f}_n(t))^2}{\sigma^2(t)}dt=0.
\end{equation*}

The diffusion coefficient $\sigma^2(\cdot)>0$ is supposed to be a known absolutely continuous function on $[0,1]$ such that 
$$\Big|\frac{d}{dt}\ln \sigma(t)\Big|\leq C,\quad t\in [0,1],$$
for some positive constant $C$.
\end{ex}

\begin{ex}\label{ex:nussbaum}
In \textcite{N96} the author establishes a global asymptotic equivalence between the problem of density estimation from an i.i.d. sample and a Gaussian white noise model. More precisely, let $(Y_i)_{i=1}^n$ be i.i.d. random variables with density $f(\cdot)$ on $[0,1]$ with respect to the Lebesgue measure. 
 The densities $f(\cdot)$ are the unknown parameters and they are supposed to belong to a certain nonparametric class $\F$ subject to a Hölder restriction: $|f(x)-f(y)|\leq C|x-y|^\alpha$ with $\alpha>\frac{1}{2}$ and a positivity restriction: $f(x)\geq \varepsilon >0$. Let us denote by $\mo_{1,n}$ the statistical model associated with the observation of the $Y_i$'s. Furthermore, let $\mo_{2,n}$ be the experiment in which one observes a stochastic process $(y_t)_{t\in[0,1]}$ such that
$$dy_t=\sqrt{f(t)}dt+\frac{1}{2\sqrt n}dW_t,\quad t\in[0,1],$$
where $(W_t)_{t\in[0,1]}$ is a standard Brownian motion.
Then the main result in \textcite{N96} is that $\Delta(\mo_{1,n},\mo_{2,n})\to 0$ as $n\to\infty$.
This is done by first showing that the result holds for certain subsets $\F_n(f_0)$ of the class $\F$ described above. Then it is shown that one can estimate the $f_0$ rapidly enough to fit the various pieces together. Without entering into any detail, let us just mention that the key steps are a Poissonization technique and the use of a functional KMT inequality.
\end{ex}

In the last years, asymptotic equivalence results have also been established for discretely observed stochastic processes. As an example, let us present the result in \textcite{esterESAIM}, very close in spirit to the one of \textcite{BL}. 
\begin{ex}Let $\{X_t\}_{t\geq0}$ be a sequence of one-dimensional time inhomogeneous jump-diffusion processes defined by
\begin{equation*}
 X_t=\eta+\int_0^tf(s)ds+\varepsilon_n\int_0^t\sigma(s)dW_s +\sum_{i=1}^{N_t}Y_i,\quad t\in [0,T],
\end{equation*}
where:
\begin{itemize}
 \item $\eta$ is some random initial condition;
 \item $W=\{W_t\}_{t\geq 0}$ is a standard Brownian motion;
 \item $N=\{N_t\}_{t\geq 0}$ is an inhomogeneous Poisson process with intensity function $\lambda(\cdot)$, independent of $W$;
 \item $(Y_i)_{i\geq 1}$ is a sequence of i.i.d. real random variables with distribution $G$, independent of $W$ and $N$;
 \item $\sigma^2(\cdot)$ is supposed to be known. The horizon of observation  $T$ is finite and $\varepsilon_n\to 0$ as $n\to\infty$.
 \item $f(\cdot)$ belongs to some non-parametric class $\F$.
 \item $\lambda(\cdot)$ and $G$ are also unknown and belong to non-parametric classes $\Lambda$ and $\mathscr G$, respectively.
\end{itemize}
In \textcite{esterESAIM}, the problem of estimating $f$ from high frequency observations of $\{X_t\}_{t\geq0}$ is considered. More precisely, we suppose to observe $\{X_t\}_{t\geq 0}$ at discrete times $0=t_0<t_1<\dots<t_n=T$ such that $\Delta_n=\max_{1\leq i\leq n}\big\{|t_{i}-t_{i-1}|\big\}\downarrow 0$ as $n$ goes to infinity. Let $\mo_n$ be the statistical model associated with the continuous observation of $\{X_t\}_{t\in [0,T]}$ and $\Qi_n$ the one associated with the observations $(X_{t_i})_{i=0}^n$. Finally, let $\Wh_n$ be the Gaussian white noise model associated with the continuous observation of the Gaussian process 
\begin{equation*}
dy_t=f(t)dt+\varepsilon_n\sigma(t)dW_t, \quad y_0 = \eta,\quad t\in[0,T].
\end{equation*}
Suppose that $\F$ is a subclass of $\alpha$-Hölder, uniformly bounded functions on $\R$ and the nuisance parameters $\sigma(\cdot)$ and $\lambda(\cdot)$ satisfy the following conditions:
\begin{itemize}
\item There exist two constants $m$ and $M$ such that $0 < m \leq \sigma(\cdot) \leq M < \infty$ and $\sigma(\cdot)$ is derivable with derivative $\sigma'(\cdot)$ in $L_{\infty}(\R)$. 
\item There exists a constant $L < \infty$ such that for all $\lambda \in \Lambda$, $\|\lambda\|_{L_2([0,T])} < L$.
\end{itemize}
Then, under the assumption that $\Delta_n^{2\alpha}\varepsilon_n^{-2}\to 0$ as $n\to\infty$, the three models $\mo_n$, $\Qi_n$ and $\Wh_n$ are asymptotically equivalent. A bound for $\Delta(\Qi_n,\Wh_n)$ and $\Delta(\mo_n,\Qi_n)$ is given by
$$\Delta_n^{\beta/2}+T\Delta_n^{2\alpha}\varepsilon_n^{-2}+T\Delta_n,$$
where $\beta=1$ if $\mathscr G$ is a subclass of discrete distributions with support on $\Z$ and $\beta=1/2$ if $\mathscr G$ is a subclass of absolutely continuous distributions with respect to the Lebesgue measure on $\R$ with uniformly bounded densities on a fixed neighborhood of $0$. In particular, this result tells us that the jumps of the process $\{X_t\}_{t\geq 0}$ can be ignored when the goal is the estimation of the drift function $f(\cdot)$.
Moreover, the proof is constructive: an explicit Markov kernel is constructed to filter the jumps out.
\end{ex}

\section{Density estimation problems and Gaussian white noise models: A constructive proof}\label{sec:density}
In this Section, following \textcite{cmultinomial} (see p. 720-725), we will detail how one can prove, in a constructive way, the asymptotic equivalence between a density estimation problem and a Gaussian white noise model, as presented in Example \ref{ex:nussbaum}. However, with respect to the work of \textcite{N96}, we will ask some stronger hypotheses on the parameter space in order to simplify the proofs. More precisely, for fixed $\gamma\in (0,1]$ and $K,\varepsilon, M$ strictly positive constants, we will consider a functional parameter space of the form 
 $$\F_{(\gamma,K,\varepsilon,M)}=\Big\{f\in C^1(I): \varepsilon\leq f(x)\leq M, \ |f'(x)-f'(y)|\leq K|x-y|^{\gamma},\ \forall x,y\in [0,1]\Big\}.$$
As in Example \ref{ex:nussbaum}, $\mo_{1,n}$ will be a density estimation problem:
\begin{equation}\label{eq:densitycarter}
(Y_i)_{1\leq i\leq n} \text{ i.i.d. r.v. with density } f\in\F_{(\gamma,K,\varepsilon,M)} 
\end{equation}
and $\mo_{2,n}$ a Gaussian with noise model:
$$dy_t=\sqrt{f(t)}dt+\frac{1}{2\sqrt n}dW_t,\quad t\in[0,1],\quad f\in\F_{(\gamma,K,\varepsilon,M)}.$$
The idea of Carter was to use the bound on the distance between multinomial and Gaussian normal variables as presented in Example \ref{ex:carter} to
make assertions about density estimation experiments. The intuition is to see the multinomial experiment as the result of grouping independent observations from a continuous density into $m$ subsets, say $J_i$, $i=1,\dots,m$. Using the square root as a variance-stabilizing transformation, these multinomial variables can be asymptotically approximated by $m$ normal variables with constant variances. These normal variables, in turn, are approximations to the increments of the process $(y_t)$ over the sets $J_i$. In Subsection \ref{sec:multinomial} we will analyze how to obtain
a asymptotically equivalent multinomial experiment starting from $\mo_{1,n}$. Assuming the results of Carter stated here as Theorems \ref{cartermultinomial} and \ref{cartergaussian} we will then obtain a bound of the $\Delta$-distance between such a multinomial experiment and one associated with 
independent Gaussian random variables. In Subsection \ref{sec:gaussian} we will explain how to show the asymptotic equivalence between an adequate normal approximation with independent coordinates and $\mo_{2,n}$.

\subsection{Density estimation problems and multinomial experiments}\label{sec:multinomial}
Let us consider a partition of $[0,1]$ in $m$ intervals $J_i=[(i-1)/m,i/m]$ and denote by $S:[0,1]^n\to \{1,\dots,n\}^m$ the application mapping the $n$-tuple $(x_1,\dots,x_n)$ to the $m$-tuple $(\#\{j:x_j\in J_1\},\dots,\#\{j:x_j\in J_m\})$, where the writing $\#\{j:x_j\in J_i\}$ stands for the number of $x_j$ belonging to the interval $J_i$. Let $P_f^{\otimes n}$ be the law of $(Y_1,\dots,Y_n)$ as in \eqref{eq:densitycarter}. The law of $S$ under $P_f^{\otimes n}$ is a multinomial distribution $\M(n;\theta_1,\dots,\theta_m)$,
$\theta_i=\int_{J_i}f(x)dx$, $i=1,\dots,m.$
In particular this means that an appropriate multinomial experiment is more informative than $\mo_{1,n}$. More precisely, we have proven that the 
statistical model associated with the multinomial distribution $(\M(n;\theta_1,\dots,\theta_m): f\in\F)$, 
denoted by $\mo_m$, is such that $\delta(\mo_{1,n},\mo_m)=0.$

Let us now investigate the quantity $\delta(\mo_m,\mo_{1,n})$.
A trivial observation is that the total variation distance between the multinomial distribution $\M(n;\theta_1,\dots,\theta_m)$ and the law $P_f^{\otimes n}$ is always $1$, hence, in order to prove that $\delta(\mo_m,\mo_{1,n})\to 0$ we need to construct a non trivial Markov kernel. We will divide the proof in three main steps.

\textbf{Step 1:} We denote by $x_i^*$ the midpoints of the intervals $J_i$, i.e. $x_i^*=\frac{2i-1}{2m}$, and we introduce a discrete random variable $X^*$ concentrated at the points $x_i^*$ with masses $\theta_i$. Let us then denote by $\mo^*$ the statistical model associated with the observation of $n$ independent realizations of $X^*$. Then, by means of a ``sufficient statistic'' argument we can get $\Delta(\mo_m,\mo^*)=0$. Indeed, consider the application
$S:\{x_1^*,\dots,x_m^*\}^n\to \{1,\dots,n\}^m$ mapping $(y_1,\dots,y_n)$ to $(\#\{j:y_j=x_1^*\},\dots,\#\{j:y_j=x_m^*\})$ and observe that the density $h$
of the law of $n$ independent realizations of $X^*$ with respect to the counting measure is given by
$$h(y_1,\dots,y_n)=\prod_{i=1}^n\p(X^*=y_i)=\theta_1^{\#\{j:y_j=x_1^*\}}\cdot\dots\cdot \theta_m^{\#\{j:y_j=x_m^*\}}.$$
By means of the Neyman-Fisher factorization theorem, we conclude that $S$ is a sufficient statistic, thus $\Delta(\mo_m,\mo^*)=0$.

\begin{figure}[!ht]
\caption{The definition of the $V_j$ functions.} \label{fig:Vj}
 \def\svgwidth{.95\textwidth}
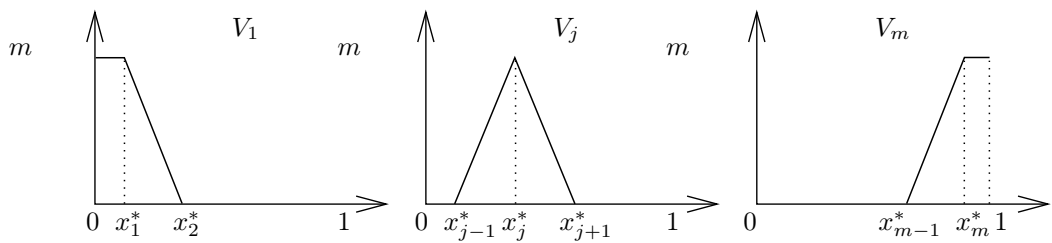
\end{figure}

\textbf{Step 2:} Starting from $n$ realizations of $X^*$ we want to obtain something close to $n$ independent realizations of $P_f$, the law of $Y_1$ as in \eqref{eq:densitycarter}. To that aim 
we define an approximation of $f$ as follows:
$$\hat f_m(x)=\sum_{j=1}^m V_j(x)\theta_j,$$
where $V_j$'s are piecewise linear functions interpolating the values in the points $x_j^*$ as in Figure \ref{fig:Vj}.

 In particular, $\hat f_m$ is a piecewise linear function that can be written as 
\begin{equation*}
 \hat f_m(x)=\begin{cases}
        m\theta_1\I_{[0,x_1^*]}(x),&\quad  \text{ if } i=1,\\ 
        \big(m-m^2|x-x_i^*|\big)\I_{[x_{i}^*,x_{i+1}^*]}(x)& \quad  \text{ if } i\in\{2,\dots,m-1\},\\
        m\theta_m\I_{[x_m^*,1]}(x),&\quad  \text{ if } i=m.
        \end{cases}
\end{equation*}
We then consider the Markov kernel 
$$M(k,A)=\sum_{j=1}^m \I_{\{x_j^*\}}(k)\int_A V_j(y)dy,\quad \forall k\in \N,\forall A\in \B([0,1]).$$
Denoting by $P^*$ the law of the random variable $X^*$, we have
\begin{align*}
 M P^*(A)=\sum_{k\in\N} M(k,A)\p(X^*=k)=\sum_{i=1}^m\theta_i M(x_i^*,A)=\sum_{i=1}^m\theta_i \int_A V_i(y)dy=\int_A \hat f_m(y)dy.
\end{align*}
Let $\hat \mo_m$ be the statistical model associated with the observation of $n$ i.i.d. random variables $(\hat Y_i)_{1\leq i\leq n}$ having
$\hat f_m $ as a density with respect to the Lebesgue measure on $[0,1]$. Applying Remark \ref{independentkernels} we get $\delta(\mo^*,\hat \mo_m)=0$. 

\textbf{Step 3:} We are only left to check that $\delta(\hat \mo_m,\mo_{1,n})\to 0$. This is actually the case and we can show that 
\begin{equation*}
 \Delta(\mo_{1,n},\hat \mo_m)=O\Big(\sqrt n\big(m^{-3/2}+m^{-1-\gamma}\big)\Big).
\end{equation*}
Indeed, the total variation distance between the family of probabilities associated with the experiments $\mo_{1,n}$ and $\hat\mo_{m}$ is bounded by
$\sqrt n H(f,\hat f_m)$.
Since $f(x)\geq \varepsilon$ for all $x\in [0,1]$ one can write:
 \begin{align*}
H^2(f,\hat f_m)&=\int_0^1\bigg(\frac{f(x)-\hat f_m(x)}{\sqrt{f(x)}+\sqrt{\hat f_m(x)}}\bigg)^2 dx\\
  &\leq \frac{1}{4\varepsilon}\int_0^1\big(f(x)-\hat f_m(x)\big)^2 dx.  
 \end{align*}
 
In order to control the $L_2$ distance between $f$ and $\hat f_m$ we will split $\int_0^1\big(f(x)-\hat f_m(x)\big)^2 dx$ as follows:
\begin{align*}
 \int_0^1\big(f(x)-\hat f_m(x)\big)^2 dx&=\int_0^{1/2m}\big(f(x)-m\theta_1\big)^2dx+\int_{1/2m}^{1-1/2m}\big(f(x)-\hat f_m(x)\big)^2dx+\\
                                        &\quad +\int_{1-1/2m}^1\big(f(x)-m\theta_m\big)^2dx.
\end{align*}
An application of the mean theorem gives $\int_{J_i}\big(f(x)-m\theta_i\big)^2dx=O(m^{-3})$, $i=1,\dots,m$. To control the term 
$\int_{1/2m}^{1-1/2m}\big(f(x)-\hat f_m(x)\big)^2dx$, let us consider the Taylor expansion of $f$ at points $x_i^*$, where $x$ denotes a point in 
$J_i$ , $i=2,\dots,m-1$:
\begin{align}\label{eq:taylor}
 f(x)&=f(x_i^*)+f'(x_i^*)(x-x_i^*)+R_i(x).
\end{align}
The smoothness condition on $f$ allows us to bound the error $R_i$ as follows: 
 \begin{align}
  |R_i(x)|&=\Big|f(x)-f(x_i^*)-f'(x_i^*)(x-x_i^*)\Big|\nonumber \\
    &=\big|f'(\xi_i)-f'(x_i^*)\big||\xi_i-x_i^*|\leq K m^{-1-\gamma}, \label{eq:resto}
 \end{align}
where $\xi_i$ is a certain point in $J_i$.

By the linear character of $\hat f_m$, we can write:
$$\hat f_m(x)=\hat f_m(x_i^*)+\hat f_m'(x_i^*)(x-x_i^*)$$
where $\hat f_m'$ denotes the left or right derivative of $\hat f_{m}$ in $x_i^*$ depending whether $x<x_i^*$ or $x>x_i^*$.
Let us observe that $\hat{f}_m'(x_i^*) = f'(\chi_i)$ for some $\chi_i \in J_i\cup J_{i+1}$ (here, we are considering right derivatives; for left ones, 
this would be $J_{i-1} \cup J_i$), indeed:
\begin{align*}
 \hat{f}_m'(x_i^*)&=-m(\hat f_m(x_i^*)-\hat f_m(x_{i+1}^*))=-m^2\bigg(\int_{\frac{i-1}{m}}^{\frac{i}{m}}f(s)ds-
 \int_{\frac{i}{m}}^{\frac{i+1}{m}}f(s)ds\bigg)\\
  &=m^2\int_{\frac{i-1}{m}}^{\frac{i}{m}}\big(f(s+1/m)-f(s)\big)ds=m\int_{J_i}f'(\xi_s)ds
\end{align*}
for some $\xi_s\in [s,s+1/m]$. Applying the mean theorem to the function $g(s)=f'(\xi_s)$ we get that $\int_{J_i}f'(\xi_s)ds=\frac{1}{m}f'(t)$ for some
$t\in J_i\cup J_{i+1}$.
The fact that $\hat{f}_m'(x_i^*)=f'(t)$, allows us to exploit the Hölder condition. Indeed, if $x\in J_i$, $i=1,\dots,m$, then there exists $t\in J_i \cup J_{i+1}$ such that: 
\begin{align*}
|f(x)-\hat f_m(x)|&\leq |f(x_i^*)-\hat f_m(x_i^*)|+|f'(x_i^*)-f'(t)||x-x_i^*|+|R_i(x)|\\
        &\leq |f(x_i^*)- \hat f_m(x_i^*)|+K|t-x_i^*||x-x_i^*|^{\gamma} + Km^{-1-\gamma}\\
        &\leq |f(x_i^*)- \hat f_m(x_i^*)|+3K m^{-1-\gamma}.
\end{align*}
Using \eqref{eq:taylor} and the fact that $\int_{J_i}(x-x_i^*)dx=0$, we get:
$$\big|f(x_i^*)-\hat f_{m}(x_i^*)\big|=m\bigg|\int_{J_i}\big(f(x_i^*)-f(x)\big)dx\bigg|\leq Km^{-1-\gamma}.$$

 Collecting all the pieces together we find
 \begin{align*}
 \int_0^1\big(f(x)-\hat f_m(x)\big)^2dx =O\big(m^{-3}+m^{-2\gamma-2}\big),
\end{align*}
hence we can conclude that $\Delta(\mo_{1,n},\hat \mo_m)=O\big(\sqrt n(m^{-3/2}+m^{-1-\gamma})\big)$.


\subsection{Independent Gaussian random variables and Gaussian white noise experiments}\label{sec:gaussian}
In Subsection \ref{sec:multinomial} we have seen how to reduce a density estimation problem to an adequate multinomial experiment. An application of the 
results of \textcite{cmultinomial} recalled in Example \ref{ex:carter} allows us to obtain an asymptotic equivalence between the statistical model 
associated with the observation of $n$ i.i.d. random variables of density $f:[0,1]\to\R$ with respect to the Lebesgue measure and an experiment 
in which one observes $m=m_n$ Gaussian and independent random variables $\Nn(\sqrt{\theta_i},1/4n)$, $i=1,\dots,m$. Of course, such a Gaussian
experiment is equivalent to $\mathscr N_m$, the statistical model associated with independent Gaussian random variables 
$\Nn\big(\sqrt{\frac{\theta_i}{m}},\frac{1}{4nm}\big)$, $i=1,\dots,m$. We claim that $\mathscr N_m$ is asymptotically equivalent to the white noise model 
$\mo_{2,n}$ associated with the continuous observation of a trajectory of a Gaussian process $(y_t)_{t\in[0,1]}$ solution of the SDE:
\begin{equation}\label{eq:wn}
dy_t=\sqrt{f(t)}dt+\frac{1}{2\sqrt n}dW_t,\quad t\in[0,1], 
\end{equation}
where $(W_t)_t$ is a standard Brownian motion. We will divide the proof in two steps. Denote by $\mathscr N_m^*$ the statistical model associated with
the observation of $(y_t)_t$ over the intervals $J_i$, $i=1,\dots,m$, i.e. $\mathscr N_m^*$ is the experiment associated with $m$ independent Gaussian
random variables $\Nn\big(\int_{J_i}\sqrt{f(s)}ds, \frac{1}{4nm}\big)$, $i=1,\dots,m$. Firstly, we will show that $\mathscr N_m$ is asymptotically equivalent
to $\mathscr N_m^*$, then that observing $(y_t)_t$ is asymptotically equivalent to observing its increments.

\textbf{Step 1:} By means of Property \ref{hellnormale} we get
\begin{align*}
 \Delta(\mathscr N_m,\mathscr N_m^*)&\leq \sqrt{2mn}\sqrt{\sum_{i=1}^m\bigg(\int_{J_i}\sqrt{f(t)}dt-\sqrt{\frac{\theta_i}{m}}\bigg)^2}\\
 &=\sqrt{2n}\sqrt{\sum_{i=1}^m \frac{1}{m}\bigg(m\int_{J_i}\big(\sqrt{f(t)}-\sqrt{m\theta_i}\big)dt\bigg)^2}\\
 \end{align*}

Denote by $E_i=|m\int_{J_i}\big(\sqrt{f(t)}-\sqrt{m\theta_i}\big)dt|$. By the triangular inequality, we bound $E_i$ by $F_i+G_i$ where:
 $$
  F_i=\bigg|\sqrt{m\theta_i}-\sqrt{f(x_i^*)}\bigg| \quad \textnormal{ and }\quad  G_i=\bigg|\sqrt{f(x_i^*)}-m\int_{J_i}\sqrt{f(y)}dy\bigg|.
 $$
 Using the same trick as in Step 3 of Subsection \ref{sec:multinomial}, we can bound:
 \begin{align*}
  F_i =\frac{|m\theta_i-f(x_i^*)|}{\sqrt{m\theta_i}+\sqrt{f(x_i^*)}}\leq \frac{|m\theta_i-f(x_i^*)|}{2\sqrt \varepsilon}=
  \frac{m}{2\sqrt \varepsilon}\Big|\int_{J_i}\big(f(s)-f(x_i^*)\big)ds\Big|=\frac{1}{2\sqrt \varepsilon}\|R_i\|_{\infty},
 \end{align*}
 where we have used the fact that $\int_{J_i}(x-x_i^*)=0$ and $R_i$ denotes the remainder in the Taylor expansion of $f$ in $x_i^*$, as in
 \eqref{eq:taylor}.
 On the other hand,
 \begin{align*}
  G_i&=m\bigg|\int_{J_j}\big(\sqrt{f(x_i^*)}-\sqrt{f(y)}\big)dy\bigg|\\
  &=m\bigg|\int_{J_i}\bigg(\frac{f'(x_i^*)}{2\sqrt{f(x_i^*)}}(x-x_i^*)+\tilde R_i(y)\bigg)dy\bigg| \leq \|\tilde R_i\|_{\infty},
 \end{align*}
 where $\tilde R_i$ is the remainder in the Taylor expansion of $\sqrt f$ in $x_i^*$. We observe that if $f$ belongs to the functional class $\F_{(\gamma,K,\varepsilon,M)}$ then $\sqrt f$ is still bounded away from zero and from infinity with a Hölder continuous derivative, more precisely $\sqrt f\in\F_{(\gamma,K/\sqrt \varepsilon,\sqrt\varepsilon,\sqrt M)}$. In particular, we deduce that $\|\tilde R_i\|_{\infty}$ has the same magnitude as $\sqrt{\frac{M}{\varepsilon}}\|R_i\|_{\infty}$. Thanks to \eqref{eq:resto} we know that $\|R_i\|_{\infty}\leq K m^{-1-\gamma}$ for any $i=2,\dots,m-1$ and
 $\|R_i\|_{\infty}=O(m^{-3/2})$ for $i\in\{1,m\}$. Hence the quantities $F_i$ and $G_i$ are of the same order and we find that
$$\Delta(\mathscr N_m,\mathscr N_m^*)=O\big(\sqrt{n}\big(m^{-1-\gamma}+m^{-\frac{3}{2}}\big)\big).$$

\textbf{Step 2:} Since $\mathscr N_m^*$ is the model associated with the observation of the increments $(\bar Y_i)_{1\leq i\leq n}$ of the process $(y_t)_t$ defined as in \eqref{eq:wn} it is clear that $\delta(\mo_{2,n},\mathscr N_m^*)=0$. Let us now discuss how to bound $\delta(\mathscr{N}_m^*,\mo_{2,n})$. We start by introducing a new stochastic process:
$$y_t^*=\sum_{i=1}^m \bar Y_i\int_0^t V_i(y)dy+\frac{1}{2\sqrt{nm}}\sum_{i=1}^mB_i(t),\quad t\in [0,1],$$
where the functions $V_i$ are defined as in Figure \ref{fig:Vj} and $B_i(t)$ are independent centered Gaussian processes independent of $(W_t)$ and with variances
$$\textnormal{Var}(B_i(t))=\int_0^tV_i(y)dy-\bigg(\int_0^tV_i(y)dy\bigg)^2.$$
These processes can be constructed from a standard Brownian bridge $B(t)$, independent of $(W_t)$, via
$$B_i(t)=B\bigg(\int_0^t V_i(y)dy\bigg).$$
By construction, $(y_t^*)$ is a Gaussian process with mean and variance given by, respectively:
\begin{align*}
 \E[y_t^*]&=\sum_{i=1}^m\E[\bar Y_i]\int_0^t V_i(y)dy=\sum_{i=1}^m\bigg(\int_{J_i}\sqrt{f(y)}dy\bigg)\int_0^t V_i(y)dy,\\
 \textnormal{Var}[y_t^*]&=\sum_{i=1}^m\textnormal{Var}[\bar Y_i]\bigg(\int_0^t V_i(y)dy\bigg)^2+\frac{1}{4nm}\sum_{i=1}^m \textnormal{Var}(B_i(t))\\
   &= \frac{1}{4nm}\int_0^t \sum_{i=1}^m V_i(y)dy= \frac{t}{4n}.
\end{align*}
One can compute in the same way the covariance of $(y_t^*)$ and deduce that 
$$Y^*_t=\int_0^t \widehat{\sqrt {f}}_m(y)dy+\int_0^t\frac{1}{2\sqrt{n}}dW^*_s,\quad t\in [0,1],$$
where
$(W_t^*)$ is a standard Brownian motion and 
$$\widehat{\sqrt {f}}_m(x):=\sum_{i=1}^m\bigg(\int_{J_i}\sqrt{f(y)}dy\bigg)V_i(x).$$

Applying Fact \ref{hellnormale}, we get that the total variation distance between the process $(y_t^*)_{t\in[0,1]}$ constructed from 
the random variables $\bar Y_i$, $i=1,\dots,m$ and the Gaussian process $(y_t)_{t\in [0,1]}$ is bounded by
$$\sqrt{4 n\int_0^1\big(\widehat{\sqrt {f}}_m(y)-\sqrt{f(y)}\big)^2dy}.$$
Since $f\in\F_{(\gamma,K,\varepsilon,M)}$ implies $\sqrt f\in \F_{(\gamma,K\sqrt M/\sqrt \varepsilon,\sqrt\varepsilon,\sqrt M)}$, the same kind of
computations made in Step 3 of Subsection \ref{sec:multinomial} allows us to conclude that
$$\Delta(\mathscr N_m^*,\mo_{2,n})=\delta(\mathscr N_m^*,\mo_{2,n})=O\big(\sqrt n\big(m^{-3/2}+m^{-1-\gamma}\big)\big).$$

\subsection{The choice of $m$}
In Subsection \ref{sec:multinomial} we have proven that the cost needed to pass from the model associated to the observation of $n$ i.i.d. random variables with unknown density $f\in\F_{(\gamma,K,\varepsilon,M)}$ to an
adequate multinomial approximation $\mathcal M(n;\theta_1,\dots,\theta_m)$ is of the order of $\sqrt n \big(m^{-3/2}+m^{-1-\gamma}\big)$. Using Theorem \ref{cartermultinomial} we can take a further step obtaining a Gaussian approximation (with independent coordinates) starting from the multinomial one.
This comes to the price of $\frac{m\ln m}{\sqrt n}$. Finally, in Subsection \ref{sec:gaussian} we have found that for appropriate choices of $m$ there is an asymptotic equivalence between such a Gaussian approximation and the Gaussian with noise model $\mo_{2,n}$. A bound for the rate of convergence of the $\Delta$-distance up to constants is, again, given by $\sqrt n \big(m^{-3/2}+m^{-1-\gamma}\big)$. 
In particular we deduce that
\begin{equation*}
\Delta(\mo_{1,n},\mo_{2,n})=\begin{cases}
                            O\Big(n^{-\frac{\gamma}{2(\gamma+2)}}\log n\Big),\quad &{if } 0<\gamma\leq\frac{1}{2},\\
                            0(n^{-\frac{1}{10}}\log n)\quad &{if } \frac{1}{2}<\gamma\leq 1.
                            \end{cases}
\end{equation*}
after the choice $m=n^{1/(2+\gamma)}$.

\printbibliography


\end{document}